\documentclass[reqno]{amsart}

\usepackage{amsmath}
\usepackage{amssymb}
\usepackage{amsfonts}
\usepackage{amsthm}
\usepackage{multirow}
\usepackage{enumerate}
 \usepackage{color}

\usepackage{array}
\newcolumntype{L}{>{\displaystyle}l}
\newcolumntype{C}{>{\displaystyle}c}
\newcolumntype{R}{>{\displaystyle}r}

\usepackage{setspace}
\newcommand{\R}{\mathbb{R}}
\newcommand{\C}{\mathbb{C}}
\newcommand{\f}{\rightarrow}
\newcommand{\deb}{\bar\partial}
\newcommand{\de}{\partial}
\newcommand{\K}{K\"{a}hler}

\newcommand{\lmb}{\lambda}
\newcommand{\ov}[1]{\overline{#1}}
\newcommand{\w}[1]{\widetilde{#1}}

\newcommand{\XX}[1]{ {\mathcal X} \left({#1}\right) }
\newcommand{\En}[2]{ \Ent\left({#1},\, {#2}\right) }

\newcommand{\D}{\mathcal{D}}

\newcommand{\hyp}{h}

\newcommand{\id}{\operatorname{id}}
\newcommand{\Id}{\operatorname{Id}}

\newcommand{\Vol}{\operatorname{Vol}}

\newcommand{\grad}{\operatorname{grad}}
\newcommand{\Ent}{\operatorname{Ent_d}}

\newcommand{\ep}{{\varepsilon}}
\newcommand{\X}{{\mathcal X}}

\newcommand{\Entv}{\operatorname{Ent_v}}
\newcommand{\di}{\rho}

\newcommand{\Entvol}{\operatorname{Ent_{vol}}}
\newcommand{\Entp}{\operatorname{Ent_{top}}}

\newtheorem{thm}{Theorem}[section]
\newtheorem{prop}{Proposition}[section]
\newtheorem{lem}[prop]{Lemma}
\newtheorem{cor}[prop]{Corollary}
\newtheorem{defn}[prop]{Definition}

\newtheorem{ex}[prop]{Example}
\newtheorem{rmk}[prop]{Remark}

\begin{document}

\date{\today}

\author[R. Mossa]{Roberto Mossa}

\address{Departamento de Matem\'atica, Universidade de S\~ao Paulo, Rua do Mat\~ao 1010, CEP 05508-900, S‹o Paulo, SP, Brazil}

\email{roberto.mossa@gmail.com}

\title[Diastatic entropy and rigidity of hyperbolic manifolds]
{Diastatic entropy and rigidity of hyperbolic manifolds}

\begin{abstract}
Let  $f: Y \f X$ be a continuous map between a compact real analytic \K\ manifold $(Y,\, g)$ and a compact complex {hyperbolic manifold} $\left( X,\, g_0 \right)$. In this paper we give a lower bound of the diastatic entropy of $(Y,\, g)$ in terms of the diastatic entropy of $\left( X,\, g_0 \right)$ and the degree of $f$. When the lower bound is attained we get geometric rigidity theorems for the diastatic entropy
analogous to the ones obtained by G. Besson, G. Courtois and S. Gallot \cite{bcg1} for the {volume entropy}. As a corollary, when $X=Y$, we show that the {minimal diastatic entropy} is achieved if and only if $g$ is isometric to the hyperbolic metric $g_0$.
\end{abstract}

\maketitle

\tableofcontents

\section{Introduction and statement of main results}
In this paper, we define the   \emph { diastatic entropy} $\Ent\left(Y,\, g \right)$ of a compact real analytic \K\ manifold $\left( Y, \, g \right)$ with \emph {globally defined diastasis function} (see Definition \ref{diast func} and \ref{defn diastatic entropy} below). This is a real analytic invariant defined, in the noncompact case, by the author in \cite{mossau}, where the link with Donaldson's balanced condition is studied. The diastatic entropy extends the concept of volume entropy using the diastasis function instead of the geodesic distance. Throughout this paper a compact \emph{complex hyperbolic manifold} will be a compact real analytic complex manifold $(X, g_0)$ endowed with locally Hermitian symmetric metric with holomorphic sectional curvature strictly negative (i.e. $\left(X, \, g_0\right)$ is the compact quotient of a complex hyperbolic space, see Example \ref{example 1} below). Our main result is the following theorem, analogous to the celebrated result of G. Besson, G. Courtois, S. Gallot on the minimal \emph{volume entropy} of a compact negatively curved locally symmetric manifold (see \eqref{eq main bcg} below) \cite[Th\'eor\`eme Principal]{bcg1}:
\begin{thm}\label{thm dentropy}
Let $\left(Y,\, g\right)$ be a compact \K\ manifold of dimension $n\geq 2$ and let $\left(X,\, g_0\right)$ be a compact complex hyperbolic manifold of the same dimension. If $f:Y \rightarrow X$ is a nonzero degree continuous map, then
\begin{equation}\label{eq main}
\Ent\left(Y,\, g \right)^{2n}\, {\Vol}\left(Y,\, g \right)\,\geq\, \left|{\deg} \left( f \right)\right| \, {\Ent}\left(X,\, g_0\right)^{2n}\, {\Vol}\left(X,\, g_0\right).
\end{equation}
Moreover the equality is attained if and only if $f$ is homotopic to a holomorphic or anti-holomorphic homothetic\footnote{$F$ is said to be homothetic if $F^* g_0 = \alpha\, g$ for some $\alpha >0$.} covering $F: Y \f X$.
\end{thm}
As a first corollary we obtain a characterization of the hyperbolic metric as that metric which realizes the minimum of the diastatic entropy:
\begin{cor}\label{cor1}
Let $\left(X,\,  g_0\right)$ be a compact complex hyperbolic manifold of dimension $n\geq 2$ and  
denote by $\mathcal E \left( X, \, g_0\right)$ the set of metrics $g$ on $X$ with globally defined diastasis and fixed volume $\Vol\left( g \right) =  \Vol\left( g_0 \right)$.
Then the functional $\mathcal F: \mathcal E \left( X, \, g_0\right) \f \R \cup \{\infty\}$ given by
$
g \stackrel{\mathcal F}{\mapsto} \Ent\left(X, \, g\right),
$
attains its minimum when $g$ is holomorphically or anti-holomorphically  {isometric} to $g_0$.
\end{cor}
This corollary can be seen as the \emph{diastatic} version of the A. Katok and M. Gromov conjecture on the minimal \emph{volume entropy} of a locally symmetric space with  strictly negative  curvature (see \cite[p. 58]{gromov}), proved by G. Besson, G. Courtois, S. Gallot in \cite{bcg1}. We  also apply Theorem \ref{thm dentropy} to give a simple proof for the complex version of the Mostow and Corlette--Siu--Thurston rigidity theorems:
\begin{cor}\label{cor2}(Mostow).
Let $\left(X, \, g_0\right)$ and  $\left(Y, \, g\right)$ be two compact complex hyperbolic manifolds of dimension $n\geq 2$. If $X$ and $Y$ are homotopically equivalent then they are holomorphically or anti-ho\-lo\-mor\-phi\-cal\-ly homothetic.
\end{cor}
\begin{cor}\label{cor3}(Corlette--Siu--Thurston).
Let $\left(X, \, g_0\right)$ and  $\left(Y, \, g\right)$ be as in the previous corollary and with the same (constant) holomorphic sectional curvature. If $f: Y\f X$ is a continuous map such that
\begin{equation}\label{eq vol}
\Vol\left( Y\right) = \left| \deg \left( f \right) \right| \Vol\left( X \right)
\end{equation}    
then there exists a holomorphically or anti-holomorphically Riemannian covering $F: Y\f X$ homotopic to $f$. 
\end{cor}
The paper consists of others two sections. In Section \ref{Diastasis, diastatic entropy and diastatic baricentre} we recall the basic definitions. Section \ref{proof cor} is dedicated to the proof of Theorem \ref{thm dentropy}. The proof is based on the analogous result for the volume entropy (see formula \eqref{eq main bcg} below) and on Lemma \ref{lower} which provides a lower bound for the diastatic entropy in terms of volume entropy. 

\smallskip

\noindent {\bf Acknowledgments.} The author would like to thank Professor Sylvestre Gallot and Professor Andrea Loi for their help and their valuable comments.

\section{Diastasis and diastasic entropy}\label{Diastasis, diastatic entropy and diastatic baricentre}

The diastasis is a special \K\ potential defined by E. Calabi in its seminal paper \cite{calabi}. Let $\left( \w Y,\, \w g \right)$ be a real analytic \K\ manifold. For every point $p \in \w Y$ there exists a real analytic function  $\Phi:V \f \R$, called \K\ potential, defined in a neighborhood $V$
of  $p$ such that
$\w  \omega =\frac{i}{2}\, \partial\ov\partial\, \Phi$, where $\w \omega$ is the \K\ form associated to $\w g$. Let $z=(z_1,\dots,z_n)$ be a local coordinates system around $p$. By duplicating the variables $z$ and $\ov z$
the real analytic \K\ potential $\Phi$ can be complex analytically
continued to a function $\hat \Phi: U\times U \f \C$  in a neighbourhood
$U\times U \subset V \times V$ of 
$\left(p, p\right)$ which is holomorphic in the first entry and antiholomorphic in the second one.
\begin{defn}[Calabi, \cite{calabi}]\label{diast func} \rm
The \emph{diastasis function} $\D:U \times U \f \R$ is defined by
\[
\D\left( z,\, w \right):=\hat \Phi \left(z,\, \ov z\right)+\hat \Phi \left(w,\,  \ov w\right)-\hat \Phi \left(z,\,  \ov w\right)-\hat \Phi \left(w,\,  \ov z\right).
\]
The \emph{diastasis function centered in $w$}, is the \K\ potential $\D_w: U \f \R$ around $w$ given by
\[
\D_w\left( z \right):=\D \left(z, \, w \right).
\]
We will say that a compact \K\ manifold $\left( Y,\, g\right)$ has \emph{globally defined diastasis} if its universal \K\ covering $\left(\w Y,\, \w g\right)$ has globally defined diastasis $\D: \w Y \times \w Y \f \R$.
\end{defn}
One can prove that the diastasis is uniquely determined by the \K\ metric $\w g$ and that it does not depend on the choice of the local coordinates system or on the choice of the \K\ potential $\Phi$.

Calabi in \cite{calabi} uses the diastasis to give necessary and sufficient conditions for the existence of an holomorphic isometric immersion of a real analytic \K\ manifolds into a complex space form. For others interesting applications of the diastasis function see \cite{LM03, LM01, LM04, LMZ01, LMZ03, LMZ02, M02} and reference therein.

Assume that $\left(\w Y,\, \w g\right)$ has globally defined diastasis $\D: \w Y \times \w Y  \f \R$. Its (normalized\footnote{Our definition of diastatic entropy differs respect to the one given in \cite{mossau} by the normalizing factor $\mathcal X \left({\w g}\right)$.}) {diastatic entropy} is defined by:
\begin{equation}\label{int entr}
\begin{split}
\Ent\left(\w Y,\, \w g\right)=\operatorname{ \mathcal X} \left({\w g}\right)\,\inf\left\{c \in \R^+: \,  \int_{\w Y} e^{-c\, \D_{w}}\, \nu_{\w g} < \infty\right\},
\end{split}
\end{equation}
where $\mathcal X \left({\w g} \right)= \sup_{y,\, z\,\in\, \w Y}{\|\grad_y \D_z\|}$ and $\nu_{\w g}$ is the volume form associated to ${\w g}$. If $\mathcal X \left({\w g}\right)=\infty$ or the infimum in \eqref{int entr} is not achieved by any $c \in \R^+$, we set $\Ent\left(\w Y,\, \w g\right)=\infty$. The definition does not depend on the base point $w$, indeed, as
\[
\left| \D_{w_1}\left(x\right) - \D_{w_2}\left(x\right) \right| = \left| \D_x\left({w_1}\right) - \D_x\left({w_2}\right) \right| \leq \X({\w g}) \, \di\left( {w_1}, \, {w_2} \right),
\]
we have 
\[
e^{-c \, \X({\w g}) \di\left( {w_1}, \, {w_2} \right)}\! \int_{\w Y} e^{-c \, \D_{w_1} \left( x \right)} \nu_{\w g}\, \leq\,  \int_{\w Y} e^{-c \, \D_{w_2} \left( x \right)} \nu_{\w g}\, \leq\, e^{c \, \X({\w g}) \di\left({w_1}, \, {w_2} \right)}\! \int_{\w Y} e^{-c \, \D_{w_1} \left( x \right)}\nu_{\w g},
\]
therefore $\int_{\w Y} e^{-c \, \D_{w_2}  \left( x \right)} \nu_{\w g} < \infty$ if and only if $\int e^{-c \, \D_{w_1} \left( x \right)}\nu_{\w g} <  \infty$.\\

\begin{defn}\label{defn diastatic entropy}\rm
Let $\left( Y,\, g\right)$ be a compact \K\ manifold with globally defined diastasis. We define the \emph{diastatic entropy} of $\left( Y, \, g \right)$ as 
\begin{equation*}
\begin{split}
\Ent\left(Y,\,  g\right)=\Ent\left( \w Y, \, \w g \right),
\end{split}
\end{equation*}
where $\left( \w Y, \, \w g \right)$ is the universal \K\ covering of $\left( Y, \, g \right)$.
\end{defn}

\begin{ex}\label{example 1}\rm
Let $\C H^n=\left\{ z \in \C^n : \|z\|^2=|z_1|^2+\dots+|z_n|^2<1 \right\}$ be the unitary disc endowed with the hyperbolic metric $\w g_{\hyp}$ of constant holomorphic sectional curvature $-4$. The associated \K\ form and the diastasis are respectively given by
\[
\w\omega_{\hyp}=-\frac{i}{2}\,\de\deb\,\log\left( 1-\|z\|^2\right) .
\] 
and
\begin{equation}\label{diast hyp}
\D^h(w,z)=-\log\left(\frac{\left(1-\|z\|^2\right)\left(1-\|w\|^2\right)}{\left|1-z w^*\right|^2}\right).
\end{equation}
Denote by $\omega_{e}=\frac{i}{2}\, \de\deb\, \| z \|^2$ the restriction to $\C H^n$ of the flat form of $\C^n$. One has
\begin{equation*}
\int_{\C H^n}e^{-\alpha \, \D^h_0}\ \frac{\omega_{\hyp}^n}{n!}
=\int_{\C H^n}\left(1-|z|^2\right)^{\alpha-n-1}\ \frac{\omega_{e}^n}{n!}<\infty\ \Leftrightarrow\ \alpha > n,
\end{equation*}
and by a straightforward computation one sees that $\XX{\w g_h}=2$. We conclude by \eqref{int entr} that  
\begin{equation}\label{diast ent hyperb}
\Ent\left( \C H^n ,\, \w g_{\hyp}\right) =2\,n. 
\end{equation}
\end{ex}
\begin{rmk}\rm
It should be interesting to compute  $\XX{g_B}$, where $g_B$ is the Bergman metric of an homogeneous bounded domain. This combined with the results obtained in \cite{mossau}, will allow us to obtain the diastatic entropy of this domains.
\end{rmk}

\section{Proof of Theorem \ref{thm dentropy} and Corollaries \ref{cor1}, \ref{cor2} and \ref{cor3}}\label{proof cor}

We start by recalling the definition of \emph{volume entropy} of a compact Riemannian manifold $(M,\, g)$. Let $\pi : \left(\w M, \, \w g\right) \f (M,\, g)$ its riemannian universal cover. We define the volume entropy of $\left( M, \, g \right)$ as
\begin{equation}\label{def vol ent int}
\Entv\left(M,\, g\right)= \inf\left\{c \in \R^+: \, \int_{\w M} e^{-c\, \w\di \left( w,\, x\right)} \, \nu_{\w g}(x)   < \infty \right\},
\end{equation}
where $\w\di$ is the geodesic distance on $\left( \w M,\, \w g\right)$ and $\nu_{\w g} $ is the volume form associated to $\w g$. By the triangular inequality, we can see that the definition does not depend on the base point $w$. As the volume entropy depends only on the Riemannian universal cover it make sense to define
\[
\Entv\left(\w M,\,\w g\right)=\Entv\left(M,\, g\right).
\]

 The \emph{classical definition} of volume entropy of  a compact riemannian manifold $\left( M, \, g \right)$, is the following
\begin{equation}\label{def vol ent class}
\Entvol \left(M,\, g\right) =  \lim_{t \f \infty} \frac{1}{t}\log\Vol\left(  B_p\left(t\right)\right),
\end{equation}
where $\Vol\left(  B_p\left(t\right)\right)$ denotes the volume of the geodesic ball $ B_p\left(t\right)\subset \w M$, of center in $p$ and radius $t$. 
This notion of entropy is related with one of the main invariant for the dynamics of the geodesic flow of $\left(M,\, g\right)$: the topological entropy $\Entp \left(M,\, g\right)$ of this flow. For every compact manifold $\left(M,\, g\right)$ A. Manning in \cite{manning} proved the inequality 
$\Entvol \left(M,\, g\right) \leq \Entp \left(M,\, g\right)$, 
which is an equality when the curvature is negative. We refer the reader to the paper \cite{bcg1} (see also \cite{bcg2} and \cite{bcg3}) of  G. Besson, G. Courtois and S. Gallot  for an overview on the volume entropy and for the proof of the celebrated minimal entropy theorem. For an explicit computation of the volume entropy $\Entv\left( \Omega,\, g\right)$ of a symmetric bounded domain $\left( \Omega,\, g\right)$ see \cite{M03}. 

The next lemma shows that the classical definition of volume entropy \eqref{def vol ent class} does not depend on the base point and it is equivalent to definition \eqref{def vol ent int}, that is
\begin{equation*}
\Entvol \left(M,\, g\right)=\Entv \left( M,\,  g\right).
\end{equation*}
\begin{lem}\label{classicalentropy}
Denote by $$\underline L := \liminf_{R \f + \infty} \left( \frac{1}{R} \log \left( \Vol B \left( x_0, \, R \right) \right) \right)$$ and $$\ov L :=\limsup_{R \f + \infty} \left( \frac{1}{R} \log \left( \Vol B \left( x_0, \, R \right) \right) \right),$$ where $B \left( x_0, \, R \right) \subset \left(\w M, \,\w g \right)$ is the geodesic ball of centre $x_0$ and radius $R$. Then the  two limits does not depends on $x_0$ and
\begin{equation*}
\underline L \leq \Entv \left( M, \, g \right) \leq  \ov L.
\end{equation*} 
\end{lem}
\proof
Let $x_1$ an arbitrary point of $M$.
Set $D = d \left( x_0, \, x_1 \right)$ and $R > D$. By the triangular inequality
\[
B \left( x_0, \, R- D \right)  \subset B \left( x_1, \, R \right) \subset B \left( x_0 ,\, R+D \right).
\] 
Let $R'=R+D$, we have
\[
\liminf_{R \f + \infty} \left( \frac{1}{R} \log \left( \Vol B \left( x_1, \, R \right) \right) \right) \leq
\liminf_{R \f + \infty} \left( \frac{1}{R} \log \left( \Vol B \left( x_0, \, R +D \right) \right) \right) 
\]
\[
= \liminf_{R' \f + \infty}  \left( \frac{R'}{R'-D}\frac{1}{R'} \log \left( \Vol B \left( x_0, \, R' \right) \right) \right) 
\]
\[
\leq\liminf_{R' \f + \infty} \left( \frac{1}{R'} \log \left( \Vol B \left( x_0, \, R' \right) \right) \right). 
\]
With the same argument one can prove the inequality in the other direction, so that $\underline L$ does not depend on $x_0$. Analogously we can prove that $\overline L$ does not depend on $x_0$.

 By the definition of limit inferior and superior, for every $\ep > 0$, there exists $R_0(\ep)$ such that, for $R \geq R_0(\ep)$,
\[
\underline L - \ep \leq \left( \frac{1}{R} \log \left( \Vol B \left( x_0, \, R \right) \right) \right) \leq \ov L + \ep
\]
equivalently
\begin{equation}\label{liminfsup}
e^{\left( \underline L - \ep \right) R}  \leq \left( \Vol B \left( x_0, \, R \right)  \right) \leq e^{\left( \ov L + \ep \right) R}.
\end{equation}
Integrating by parts we obtain
\[
I:=\int_{\w M} e^{-c\,  \w \di \left( x_0 ,\, x \right)} dv(x) = \int_0^\infty e^{-c\, r} \Vol_{n-1} \left( S\left( x_0 ,\, r \right) \right) dr 
\]
\[
= \Vol \left( B\left( x_0 ,\, r \right) \right) e^{-c\,r} \Big |_0^{\infty} + c \int_0^\infty e^{-c\,  r}  \Vol \left( B\left( x_0 ,\, r \right) \right) dr.
\]
where $S\left( x_0 ,\, r \right)=\de B\left( x_0 ,\, r \right)$. On the other hand, by \eqref{liminfsup} we get
\[
\int_{R_0(\ep)}^\infty e^{\left(  \underline L  -c -\ep \right)  r}\,  dr \leq \int_{R_0(\ep)}^\infty e^{-c\,  r}  \Vol \left( B\left( x_0 ,\, r \right) \right) dr \leq \int_{R_0(\ep)}^\infty e^{-\left( c - \ov L -\ep   \right)  r}\,  dr.
\]
We deduce that if $c > \ov L$ then $I$ is convergent i.e $\ov L \geq \Entv$ and that if $I$ is not convergent when $c < \underline L$, that is $\Entv \geq \underline L$, as wished. 
\endproof

The next lemma show that the diastatic entropy is bounded from below by the volume entropy. 
\begin{lem}\label{lower} Let $\left(Y,\, g\right)$ be a compact \K\ manifold with globally defined diastasis, then
\begin{equation}\label{eq thm upper}
\Ent\left(Y,\, g \right)\,\geq\,\Entv\left(Y,\, g \right).
\end{equation}
This bound is sharp when $\left(Y,\, g\right)$ is a compact quotient of the complex hyperbolic space.  That is,
\begin{equation}\label{eq diast vol ent}
\Ent\left(\C H^n,\, \w g_h \right)\,=2\,n=\,\Entv\left(\C H^n,\, \w g_h \right).
\end{equation}
\end{lem}
\proof 
Let $(\w Y,\,\w g)$ be universal \K\ cover of $(Y,\, g)$. For every $w,x \in \w Y$ we have 
\[
\D_w\left(x\right) = \D_w\left( x\right) - \D_w\left(w\right) \leq \sup_{z \in \w Y} \left\| d_ z \D_w \right\| \di_w \left( x \right) \leq {\X ({\w g})}\ \di_w\left( x \right),
\]
so
\[
\int_{\w Y} e^{-c\, \X \left(\w g \right) \, \di_w\left( x \right)}\, \nu_{\w g}  \leq \int_{\w Y} e^{-c\, \D_w \left( x \right)}\, \nu_{\w g}.
\]
Therefore, if $c \, {\X ( {\w g} )} \leq \Entv\, ({\w Y},\, {\w g}) $ then $c\, {\X \,( {\w g} )} \leq \Ent \, ({\w Y},\, {\w g} )$. We obtain \eqref{eq thm upper} by setting $c=\frac {\Entv \, ({\w Y},\, {\w g})}{{ \X \left(\w g \right)}}$. Equation \eqref{eq diast vol ent} follow by \eqref{diast ent hyperb} and \cite[Theorem 1.1]{M03}.
\endproof

\noindent
{\bf Proof of Theorem \ref{thm dentropy}.}
Let $\left(X,\, g_0 \right)$ as in Theorem \ref{thm dentropy} and let $\pi_X : \left(\C H^n,\, \w g_0 \right)\f \left(X,\, g_0 \right)$ be the universal covering. Notice that $\w g_0 = \lmb \, \w g_h$ for some positive $\lmb$. Then we have
\begin{equation}\label{eq entd entv}
\begin{array} {C}
{\Vol}\left(X,\, g_0\right)\Entv\left(X,\,g_0\right)^{2n}={\Vol}\left(X,\, g_h\right)\Entv\left(X,\,g_h\right)^{2n} \\[1em]
={\Vol}\left(X,\, g_h\right)\Ent\left(X,\, { g_h}\right)^{2n}={\Vol}\left(X,\, g_0\right)\Ent\left(X,\, { g_0}\right)^{2n},
\end{array}
\end{equation}
where the first and the third equality are consequence of the fact that $\Entv\left(\C H^n,\,\w g_0\right)= \frac {1}{\sqrt \lmb} \Entv\left(\C H^n,\,\w g_h\right)$ and $\Ent\left(\C H^n,\,\w g_0\right)= \frac {1}{\sqrt \lmb} \Ent\left(\C H^n,\,\w g_h\right)$, while the second equality follows by \eqref{eq diast vol ent}. Let $f: Y \f X$ be as in Theorem \ref{thm dentropy}, then, by  \cite[Th\'eor\`eme Principal]{bcg1} we know that
\begin{equation}\label{eq main bcg}
\Entv\left(Y,\, g \right)^{2n}\, {\Vol}\left(Y,\, g \right)\,\geq\, \left|{\deg} \left( f \right)\right| \, {\Entv}\left(X,\, g_0\right)^{2n}\, {\Vol}\left(X,\, g_0\right)
\end{equation}
where the equality is attained if and only if $f$ is homotopic to a homothetic covering $F: Y \f X$. Putting together   \eqref{eq thm upper}, \eqref{eq entd entv} and \eqref{eq main bcg} we get that
\begin{equation*}
\Ent\left(Y,\, g \right)^{2n}\, {\Vol}\left(Y,\, g \right)\,\geq\, \left|{\deg} \left( f \right)\right| \, {\Ent}\left(X,\, g_0\right)^{2n}\, {\Vol}\left(X,\, g_0\right)
\end{equation*}
where the equality is attained if and only if $f$ is homotopic to a homothetic covering $F: Y \f X$.

To conclude the proof it remains to prove that $F$ is holomorphic or anti-holo\-morph\-ic.
Up to homotheties, it is not restrictive to assume that $g = F^*{ g_0}$, so that its lift $\w F:\w Y \f  \C H^n$ to the universal covering it is a global isometry. Fix a point $q \in \w Y$, let $p = {\w F}(q)$ and denote $A_q=\w F^*{J_0}_p$ the endomorphism acting on $T_{q}\w Y$, where $J_0$ is the complex structure of $\C H^n$. Denote by $\mathcal G_{\w Y}$ and respectively $\mathcal G_{\C H^n}$ the holonomy groups of $( \w Y,\, \w g )$ and respectively $\left( \C H^n,\, \w g_0\right)$. Note that $\mathcal G_{\w Y}={\w F}^*\mathcal G_{\C H^n}$ and that $\mathcal G_{\C H^n}=SU(n)$, therefore $\mathcal G_{\w Y}$ acts irreducibly on $T_q \w Y$. As $J_0$ commutes with the action of $\mathcal G_{\C H^n}$, by construction $A_q$ is invariant with respect to the action of $\mathcal G_{\w Y}$. Therefore, denoted $\Id_q$ the identity map of $T_q \w Y$, by  Schur's lemma, $A_q=\lmb \Id_q$ with $\lmb \in \C$. Moreover $-\Id_q=A_q^2=\lmb^2 \Id_q$, so $\lmb=\pm i$. By the arbitrarity of $q$ we conclude that $\w F$ is holomorphic or anti-holomorphic.

\medskip

\noindent {\bf Proof of Corollary \ref{cor1}.} This is an immediate consequence of Theorem \ref{thm dentropy} once assumed $Y = X$, $\Vol \left( g \right) = \Vol \left( g_0 \right)$ and $f=\id_X$ the identity map of $X$.

\medskip

\noindent {\bf Proof of Corollary \ref{cor2}.} Let $h: Y \f X$ be an homotopic equivalence and $h^{-1}$ its homotopic inverse. Substituting in \eqref{eq main}, once with $f=h$ and once with $f=h^{-1}$, we have respectively
\[
\Ent\left(Y,\, g \right)^{2n}\, {\Vol}\left(Y,\, g \right)\,\geq \, \left|{\deg} \left( h \right)\right| \,  {\Ent}\left(X,\, g_0\right)^{2n}\, {\Vol}\left(X,\, g_0\right)
\]
and
\[
  {\Ent}\left(X,\, g_0\right)^{2n}\, {\Vol}\left(X,\, g_0\right)\,\geq \, \left|{\deg} \left( h^{-1} \right)\right| \,\Ent\left(Y,\, g \right)^{2n}\, {\Vol}\left(Y,\, g \right).
\]
We then conclude that  $\Ent\left(Y,\, g \right)^{2n}\, {\Vol}\left(Y,\, g \right)\,= \,{\Ent}\left(X,\, g_0\right)^{2n}\, {\Vol}\left(X,\, g_0\right)$ and that $\left|{\deg} \left( h \right)\right|=1$. Therefore, by applying the last part of Theorem \ref{thm dentropy}, we see that $h$ is homotopic to an holomorphic (or antiholomorphic)  homothety $F:X \f Y$.

\medskip

\noindent {\bf Proof of Corollary \ref{cor3}.} Let $\pi_Y : \left(\C H^n,\, \w g \right)\f \left(Y,\, g \right)$ and $\pi_X : \left(\C H^n,\, \w g_0 \right)\f \left(X,\, g_0 \right)$ be the universal coverings, since $g_0$ and $g$ are both hyperbolic with the same curvature, we conclude that $\w g_0= \w g$ and that $\En{X}{g_0}=\En{Y}{g}$. Therefore we get an equality in \eqref{eq main}. Using again the last part of Theorem \ref{thm dentropy} we get $\Vol\left( Y\right) = \left| \deg \left( F \right) \right| \Vol\left( X \right)$ and we conclude that $F$ is locally isometric.

\end{document}